\documentclass[12pt]{amsart}

\usepackage[latin1]{inputenc}
\usepackage{amsmath}
\usepackage{amsfonts}
\usepackage{amssymb}
\usepackage{graphics}
\usepackage{amssymb,amsmath,amsthm,amscd,epsf,latexsym,verbatim,graphicx,amsfonts}
\input epsf.tex

\newtheorem{theorem}{Theorem}[section]

\numberwithin{equation}{section}

\newcommand{\bbR}{\mathbb{R}}
\newcommand{\bbC}{\mathbb{C}}

\newcommand{\bbD}{\mathbb{D}}
\newcommand{\bbN}{\mathbb{N}}

\newcommand{\dtpi}{\frac{d\theta}{2\pi}}

\newcommand{\eitheta}{e^{i\theta}}

\newcommand{\mcn}{\mathcal{N}}

\newcommand{\bard}{\overline{\bbD}}

\newcommand{\supp}{\textrm{supp}}

\newcommand{\nri}{n\rightarrow\infty}

\begin{document}

\title[ ] {Weak Convergence of CD Kernels: A New Approach on the Circle and Real Line}

\bibliographystyle{plain}

\thanks{  }

\vspace{-7mm}

\maketitle

\begin{center}
\textbf{Brian Simanek}\footnote{Mathematics MC 253-37, California Institute of Technology, Pasadena, CA 91125, USA. E-mail: bsimanek@caltech.edu.  Supported in part by an NSF GRFP grant.}
\end{center}

\begin{abstract}
Given a probability measure $\mu$ supported on some compact set $K\subseteq\bbC$ and with orthonormal polynomials $\{p_n(z)\}_{n\in\bbN}$, define the measures
\[
d\mu_n(z)=\frac{1}{n+1}\sum_{j=0}^n|p_j(z)|^2d\mu(z)
\]
and let $\nu_n$ be the normalized zero counting measure for the polynomial $p_n$.  If $\mu$ is supported on a compact subset of the real line or on the unit circle, we provide a new proof of a 2009 theorem due to Simon that for any fixed $k\in\bbN$ the $k^{th}$ moment of $\nu_{n+1}$ and $\mu_n$ differ by at most $O(n^{-1})$ as $\nri$.
\end{abstract}

\vspace{4mm}


\vspace{2mm}


\vspace{2mm}

\normalsize

\section{Introduction}\label{intro}

Given a probability measure $\mu$ on $\bbC$ with infinite and compact support, we can form the sequence $\{p_n(z)\}_{n\in\bbN}$ of orthonormal polynomials satisfying
\[
\int_{\bbC}\overline{p_n(z)}p_m(z)d\mu(z)=\delta_{n,m}
\]
and normalized so that each $p_n$ has positive leading coefficient $\kappa_n$.  With this sequence, we define
\[
K_n(z,\zeta;\mu)=\sum_{j=0}^n\overline{p_j(\zeta)}p_j(z),
\]
the so-called \textit{Reproducing Kernel} for polynomials of degree $n$.  We assign it this name because of the reproducing property, namely that if $Q$ is any polynomial of degree at most $n$ then
\[
Q(w)=\int_{\bbC}Q(z)K_n(w,z;\mu)d\mu(z).
\]
With this notation, we can define the probability measures
\[
d\mu_n=\frac{K_{n}(z,z;\mu)}{n+1}d\mu
\]
for each $n\in\bbN$.

If we write $p_n(z)=\kappa_n\prod_{j=1}^n(z-z_j^{(n)})$ (the $z_j^{(n)}$ need not be distinct), then we define the measures
\[
d\nu_n=\frac{1}{n}\sum_{j=1}^n\delta_{z_j^{(n)}}
\]
for each $n\in\bbN$.  In \cite{WeakCD}, Simon proved the following theorem:

\begin{theorem}\label{simon}\cite{WeakCD}
Let $N(\mu)=\sup\{|z|:z\in\supp(\mu)\}$.  For any $k\in\bbN$, we have
\begin{align}\label{main}
\left|\int_{\bbC}z^{k}d\mu_n(z)-\int_{\bbC}z^kd\nu_{n+1}(z)\right|\leq\frac{2kN(\mu)^{k}}{n+1}.
\end{align}
\end{theorem}

From this theorem, Simon deduces the following important corollary.  Suppose $C$ is a circle centered at $0$ with radius larger than $N(\mu)$.  Let $\hat{\mu}_n$ denote the balayage (see Theorem II.4.1 in \cite{SaffTot}) of the measure $\mu_n$ onto $C$ and similarly define $\hat{\nu}_n$.  It follows from Theorem \ref{simon} that
\[
w\textrm{-}\lim_{\nri \atop n\in\mcn}\hat{\nu}_{n+1}=\sigma\quad\iff\quad
w\textrm{-}\lim_{\nri \atop n\in\mcn}\hat{\mu}_{n}=\sigma,
\]
where $\sigma$ is a probability measure on $C$ and $\mcn\subseteq\bbN$ is a subsequence.

The proof of the above theorem in \cite{WeakCD} relies on the relationship between the polynomials $\{p_n\}_{n\in\bbN}$ and the eigenvalues of the operator $M_z$ acting on $L^2(\mu)$ by $M_z(f(z))=zf(z)$.  In this paper, we will provide a new proof of this theorem when $\supp(\mu)\subseteq\bbR$ or $\supp(\mu)\subseteq\partial\bbD:=\{z:|z|=1\}$.  The key idea will be to look at Pr\"{u}fer phases of the appropriate ratio of the orthonormal polynomials.

If $\mu$ is supported on $\partial\bbD$, we define $\eta_n(\theta):[0,2\pi]\rightarrow\bbR$ to be a continuous function so that
\begin{align}\label{blash}
e^{i\eta_n(\theta)}=\frac{p_{n+1}(\eitheta)}{p_{n+1}^*(\eitheta)},
\end{align}
where $p_{n+1}^*(z)=z^{n+1}\overline{p_{n+1}(\bar{z}^{-1})}$ (so that the right hand side of (\ref{blash}) is a Blaschke product).  If $\mu$ is supported on $\bbR$ (we always assume compact support), then we may define $\theta_n(x):\bbR\rightarrow(-\pi/2,\infty)$ to be a continuous function so that
\begin{align}\label{realtan}
\tan(\theta_n(x))=\frac{a_np_n(x)}{p_{n-1}(x)},
\end{align}
(see Proposition 6.1 in \cite{FineIV}) where $a_n$ is a positive real number so that $p_{n-1}$ and $a_np_n$ have the same leading coefficient.  In our proofs, we will use the functions $\eta_n$ and $\theta_n$ (more precisely their derivatives) to obtain measures that approximate the measure $\mu$ in a sense suitable for our purposes.

More precisely, two approximating measures will enter.  In the unit circle case, we define
\begin{align}\label{BZ}
d\mu^n=|p_{n+1}(\eitheta)|^{-2}\dtpi
\end{align}
for each $n\in\bbN$.
The measure $\mu^n$ (called the $n^{th}$ \textit{Bernstein-Szeg\H{o}} measure) is in fact a probability measure on $[0,2\pi]$ and it induces a measure on $\partial\bbD$ with the same first $n$ moments - and hence the same first $n$ orthonormal polynomials - as $\mu$ (this follows from Theorems 1.7.8 and 1.5.5 in \cite{OPUC1}).  In the real line case, we define
\begin{align}\label{rhon}
d\rho_n=\frac{dx}{\pi(a_{n+1}^2p_{n+1}(x)^2+p_n(x)^2)}
\end{align}
as in Theorem 2.1 in \cite{CRM}.  It follows from equation (2.7) in \cite{CRM} that $d\rho_n$ is a probability measure and
\begin{align}\label{realmom}
\int_{\bbR}x^{\ell}d\rho_n(x)=\int_{\bbR}x^{\ell}d\mu(x)\quad,\quad\ell=0,1,\ldots,2n.
\end{align}

In the next section we provide our new proof of Theorem \ref{simon} when $\mu$ is supported on the unit circle.  In Section \ref{real} we consider $\mu$ supported on the real line and prove Theorem \ref{simon} with the right hand side of (\ref{main}) replaced by $O(n^{-1})$.

\section{The Unit Circle Case}\label{circle}

Our goal in this section is to provide a new proof of Theorem \ref{simon} when $\mu$ is supported on the unit circle.  We begin our proof by noting that the theorem is equivalent to the statement that the moments of the signed measures $d\hat{\nu}_{n+1}-d\mu_n$ converge $0$ at a certain rate where $\hat{\nu}_n$ is the balayage of the measure $\nu_n$ onto $\partial\bbD$.  It is easy to check that (see equation (8.2.8) in \cite{OPUC1})
\[
d\hat{\nu}_{n+1}=\frac{1}{n+1}\sum_{j=1}^{n+1}\frac{1-|z_j^{(n+1)}|^2}{|\eitheta-z_j^{(n+1)}|^2}\dtpi.
\]
If we define $\eta_n:[0,2\pi]\rightarrow\bbR$ as in (\ref{blash}) above, then equation (6.10) in \cite{Stoiciu} implies that
\[
\frac{d}{d\theta}\eta_n(\theta)=\sum_{j=1}^{n+1}\frac{1-|z_j^{(n+1)}|^2}{|\eitheta-z_j^{(n+1)}|^2}.
\]
Furthermore, equation (10.8) in \cite{FineIV} tells us that
\[
\frac{d}{d\theta}\eta_n(\theta)=\frac{K_n(\eitheta,\eitheta;\mu)}{|p_{n+1}(\eitheta)|^2}
\]
so we conclude that
\[
d\hat{\nu}_{n+1}=\frac{K_n(\eitheta,\eitheta;\mu)}{n+1}d\mu^n(\theta).
\]
Therefore, if $k\in\bbN$, we can write
\begin{align*}
\int_{\bard}z^k\,d\nu_{n+1}(z)-\int_{\partial\bbD}z^kd\mu_n(z)&=
\frac{1}{n+1}\sum_{j=0}^{n}\left[\langle p_j(z),z^kp_j(z)\rangle_{\mu^n}-
\langle p_j(z),z^kp_j(z)\rangle_{\mu}\right].
\end{align*}
Since $\mu$ and $\mu^n$ have the same first $n$ moments, at most $k$ of these summands are non-zero and each non-zero summand has absolute value at most $2$.  We have therefore proven
\[
\left|\int_{\bard}z^k\,d\mu_n(z)-\int_{\partial\bbD}z^k\,d\nu_{n+1}(z)\right|\leq\frac{2k}{n+1}
\]
exactly as in Theorem \ref{simon}.

\vspace{2mm}

\noindent\textbf{Example.}  Let $\mu$ be the normalized arclength measure on the unit circle.  In this case we have $p_n(z)=z^n$ for all $n$ and $\mu_n=\mu$ for all $n$.  The measures $\nu_n$ are all simply the point mass at $0$ with weight $1$.  This example illustrates the fact that in general, the measures $\mu_n$ and $\nu_n$ need not resemble each other as measures on $\overline{\bbD}$, so it really is important that we consider the balayage.

\section{The Real Line Case}\label{real}

Our goal in this section is to provide a new proof of Theorem \ref{simon} when $\mu$ is supported on a compact subset of the real line and with the right hand side of (\ref{main}) replaced by $O(n^{-1})$ where the implied constant depends on $k$.  There is a proof of this result due to Totik, also appearing in \cite{WeakCD}, but with the right hand side of (\ref{main}) replaced by $o(1)$ (though it can be modified to produce the same $O(n^{-1})$ discrepancy estimate for the moments as in (\ref{main}) above).
Totik's proof uses Gaussian quadratures and the monotonicity (in $n$) of the sequence $K_n(x,x;\mu)$ to establish the weak convergence result for all polynomials that are positive on the convex hull of the support of $\mu$.  The proof we present here will be analogous to the proof in Section \ref{circle} and will rely on the sequence of approximating measures $\rho_n$ (see (\ref{rhon}) above).  We will make use of formula (\ref{realrel}) below, which relates a set of perturbed zero-counting measures to a set of perturbed quadrature measures.  Combining this with an interlacing property will allow us to derive the $O(n^{-1})$ estimate in (\ref{main}).

Our computation will be a bit longer than in the unit circle case partly because in Section \ref{circle}, the most difficult calculation was already done for us in \cite{FineIV} and partly because the high moments of the measure $\rho_n$ defined in equation (\ref{rhon}) are infinite, so we need to use a cutoff function.

Let us assume $\mu$ has support contained in $[-M,M]$ and define
\[
\tau(x)=\chi_{[-M-1,M+1]}(x).
\]
Corresponding to $\mu$ there is a Jacobi matrix $J$, which is the matrix of multiplication by $x$ in the Hilbert space $L^2(\mu)$ with respect to the basis given by the orthonormal polynomials.  For any $\lambda\in\bbR$, we will let $\mu_{n,\lambda}$ be the spectral measure corresponding to the Jacobi matrix $J_n+\lambda\langle e_n,\cdot\rangle e_n$ and the vector $e_1$ where $J_n$ is the upper left $n\times n$ block of $J$ (see Section 6 in \cite{Darboux}).  Notice that $\mu_{n,\lambda}$ is supported on $n$ distinct points and by Corollary 6.3 in \cite{Darboux}, the points in the support of $\mu_{n,\lambda}$ interlace for distinct values of $\lambda$.  Let $\nu_{n,\lambda}$ be the measure placing weight $n^{-1}$ on each point in the support of $\mu_{n,\lambda}$ (so that $\nu_{n,0}=\nu_n$).  It follows from formula (6.16) in \cite{Darboux} that
\begin{align}\label{realrel}
\frac{1}{n}d\mu_{n,\lambda}=\frac{1}{K_{n-1}(x,x;\mu)}d\nu_{n,\lambda}.
\end{align}
Therefore for any fixed $k\in\bbN$, we have
\begin{align}\label{start}
\int_{\bbR}x^k\tau(x)d\nu_{n+1,\lambda}=\frac{1}{n+1}\int_{\bbR}x^k\tau(x)K_{n}(x,x;\mu)d\mu_{n+1,\lambda}.
\end{align}
After taking a suitable average (in $\lambda$), the expression on the left-hand side of (\ref{start}) approximates the $k^{th}$ moment of $\nu_{n+1}$ as $\nri$ while the right-hand side approximates the $k^{th}$ moment of $\mu_n$ as $\nri$.
Indeed, our first step is to integrate the left hand side of (\ref{start}) from $-\infty$ to $\infty$ with respect to $\frac{d\lambda}{\pi(1+\lambda^2)}$.  Notice that for any value of $\lambda$, at most one point in the support of $\nu_{n+1,\lambda}$ lies outside $[-M-1,M+1]$ because of the interlacing property.  Therefore, we have
\begin{align}\label{leftest}
\int_{-\infty}^{\infty}\int_{\bbR}x^k\tau(x)d\nu_{n+1,\lambda}(x)\frac{d\lambda}{\pi(1+\lambda^2)}=
\int_{\bbR}x^k\tau(x)d\nu_{n+1,0}(x)+O(n^{-1})
\end{align}
as $\nri$.

If we integrate the right hand side of (\ref{start}) in the same way, this becomes
\begin{align}\label{tworight}
\frac{1}{n+1}\int_{\bbR}x^k\tau(x)K_{n}(x,x;\mu)\,d\rho_n(x)
\end{align}
by Theorem 2.1 in \cite{CRM}.  Notice that this integral would be infinite without the cut-off function $\tau$.  As an aside, we note that by Proposition 6.1 in \cite{FineIV}, (\ref{tworight}) is just
\[
\frac{1}{n+1}\int_{\bbR}x^k\tau(x)\frac{1}{\pi}\frac{d\theta_{n+1}(x)}{dx}\,dx,
\]
which is why we call this the analog of the proof in Section \ref{circle}.
Notice that for any fixed $m\leq n$ we have
\[
\int_{\bbR}x^k\tau(x)|p_m(x)|^2\,d\rho_n(x)\leq(M+1)^k.
\]
This follows from the fact that $p_m$ is also the degree $m$ orthonormal polynomial for the measure $d\rho_n$ by (\ref{realmom}).  Therefore, we can rewrite (\ref{tworight}) as
\begin{align*}
\frac{1}{n+1}\int_{\bbR}x^k\tau(x)K_{n-k}(x,x;\mu)\,d\rho_n(x)+O(n^{-1})
\end{align*}
as $\nri$.
We can rewrite this again as
\begin{align}\label{threeright}
\frac{1}{n+1}\int_{\bbR}x^kK_{n-k}(x,x;\mu)\,d\rho_n(x)-\frac{1}{n+1}\int_{|x|>M+1}x^kK_{n-k}(x,x;\mu)\,d\rho_n(x)
+O(n^{-1})
\end{align}
as $\nri$.  Notice that $x^kK_{n-k}(x,x;\mu)$ is a polynomial of degree $2n-k$ while the denominator of the weight defining the measure $\rho_n$ is a polynomial of degree $2n+2$.  Therefore, both integrals in (\ref{threeright}) are finite.  The first term in (\ref{threeright}) is equal to
\[
\frac{1}{n+1}\int_{\bbR}x^kK_{n}(x,x;\mu)d\mu(x)+O(n^{-1})=
\int_{\bbR}x^k\,d\mu_n(x)+O(n^{-1})
\]
as $\nri$ again by (\ref{realmom}).  We will be finished if we can show that the second term in (\ref{threeright}) tends to $0$ like $O(n^{-1})$ as $\nri$ and for this it suffices to put a uniform bound on
\begin{align}\label{bound}
\int_{|x|>M+1}x^kK_{n-k}(x,x;\mu)\,d\rho_n(x).
\end{align}

To do this, we rewrite (\ref{bound}) as
\[
\int_{-\infty}^{\infty}\int_{|x|>M+1}x^kK_{n-k}(x,x;\mu)\,d\mu_{n+1,\lambda}(x)\,\frac{d\lambda}{\pi(1+\lambda^2)}.
\]
Recall that for each fixed $\lambda$, at most one point in the support of $\mu_{n+1,\lambda}$ has absolute value larger than $M+1$.  Let us denote this point (if it exists) by $x_{n+1,\lambda}$.  Therefore, the above integral is just
\begin{align}\label{rightfour}
\int_{A}x_{n+1,\lambda}^k\frac{K_{n-k}(x_{n+1,\lambda},x_{n+1,\lambda};\mu)}
{K_{n}(x_{n+1,\lambda},x_{n+1,\lambda};\mu)}\frac{d\lambda}{\pi(1+\lambda^2)},
\end{align}
where we used (\ref{realrel}) and the integral is taken over some set $A\subseteq\bbR$ such that $x_{n+1,\lambda}$ exists if and only if $\lambda\in A$.  Using the Christoffel Variational Principle (Theorem 9.2 in \cite{Darboux}), it is easily seen that
\[
\frac{K_{n-k}(x_{n+1,\lambda},x_{n+1,\lambda};\mu)}
{K_{n}(x_{n+1,\lambda},x_{n+1,\lambda};\mu)}\leq\left(\frac{M}{|x_{n+1,\lambda}|}\right)^{2k}.
\]
Therefore, we can bound (\ref{rightfour}) from above in absolute value by
\[
\int_A\frac{M^{2k}}{|x_{n+1,\lambda}|^k}\frac{d\lambda}{\pi(1+\lambda^2)},
\]
which is uniform in $n$ since $|x_{n+1,\lambda}|>M+1$.  This completes the proof.

\vspace{4mm}

\noindent\textbf{Acknowledgments.}  It is a pleasure to thank my advisor Barry Simon for many useful comments and suggestions.


\end{document}